\magnification= \magstep1
\input coman.tex
\input amssym.def

\riferimentifuturi
\indice

\autobibliografia

\biblitem{Coulhon-Grigoryan}
       T.~Coulhon, A.~Grigoryan, {\it Random walks on graphs with regular
	volume growth}, GAFA {\bf 8} (1998), 656-701.

\biblitem{Zucca1}
       F.~Zucca, {\it The mean value property for harmonic functions on graphs and trees}, appearing on Ann.~Mat.~Pura.~App.

\biblitem{Hille}
        E.~Hille, {\it Analytic Function Theory: vol.~I}, Chelsea Publ.~co., 
			New York, N.Y.~(1959).

\biblitem{Rudin2}
        W.~Rudin, {\it Principles of Mathematical Analysis}, Mc Graw-Hill,  (1953).

\biblitem{Kanai1} M.~Kanai, {\it Rough isometries and combinatorial
			approximations of geometries of non-compact Riemannian
			manifolds.}, J.~Math.~Soc.~Japan {\bf 37} (1985), 
			391-413 (4, 6, 9).

\biblitem{Kanai2} M.~Kanai, {\it Rough isometries and the parabolicity
			of Riemannian manifolds.}, J.~Math. Soc.~Japan {\bf 38} 			(1986), 227-238 (1, 4, 9).

\biblitem{Mark-Mcgu-Thom} St.~Markvorsen, S.~McGuinness, C.~Thomassen, 
			{\it Transient random walks on graphs and metric 
			spaces, with applications on hyperbolic surfaces.},
			Proc.~London.~Math.~Soc.~{\bf 64} (1992), 1-20 (4, 9).

\biblitem{Lyon-Sull} T.~Lyons, D.~Sullivan, {\it Function theory, 
			random path and covering spaces.}, 
			J.~Diff. Geom.~{\bf 19}	(1984), 299-323 (4, 7, 9).

\biblitem{Lyons} T.~Lyons, {\it Instabilty of the Liouville property
			for quasi-isometric Riemannian manifolds and 
			reversible Markov Chains.}, J.~Diff.~Geom.~{\bf 26} 	
			(1987), 33-66 (8, 9).

\biblitem{Berta1} D.~Bertacchi, {\it Probabilistic properties of 
			$DL$-graphs.}, preprint.

\biblitem{Fig-Picar} A.~Fig\`a-Talamanca, M.~A.~Picardello, {\it Harmonic analysis on
	free groups}, Lecture Notes in Pure and Appl.~Math., vol.~87, Dekker, New York
	and Basel, 1987.

\biblitem{Cartier2} P.~Cartier, {\it Harmonic analysis on trees}, 
	Proc.~Sympos.~Pure~~Math.~, vol.~26, Amer. Math.~Soc., Providence,
	R.I.,  1972, 419-424.

\biblitem{ErFePo} P.~Erd\"os, W.~Feller, H.~Pollard,
        {\it A theorem on power series}, Bull.~Amer.~Soc.~{\bf 55},
         (1949) 201-204.

\biblitem{Cartier} P.~Cartier, {\it Fonctions harmoniques sur un arbre}, Symposia
        Math.~, {\bf 9} (1972) 203-270.

\biblitem{DB&FZ2}
        D.~Bertacchi, F.~Zucca, {\it Asymptotic uniform estimates for
        transition probabilities on 2-comb},
        preprint.

\biblitem{Woess3}     
        W.~Woess, {\it Catene di Markov e Teoria del Potenziale nel
        Discreto}, Quaderno U.M.I. {\bf 41}, Ed.~Pitagora, Bologna (1996).

\biblitem{Engelking}     
        R.~Engelking, {\it General Topology}, Sigma series in pure 	mathematics {\bf 6}, Heldermann, Berlin (1989).

\biblitem{Pic-Woess2}     
        M.~Picardello, W.~Woess, {\it A converse to the mean value property
        on homogeneous trees}, Trans. of A.M.S.~ {\bf 311} (1999), No.~1,
        209-225.

\biblitem{Doob1}
        J.~L.~Doob, {\it Discrete potential theory and boundaries}, 
J.~Math.~Mech.~{\bf 8}, (1959), 433-458.

\biblitem{Dynkin1}
	E.~B.~Dynkin, {\it The boundary theory of Markov processes
(discrete case)}, Uspehi Mat. Nauk {\bf 24} (1969) no. 2 (146) 3--42.

\biblitem{Singer1}
        I.~Singer, {\it Bases in Banach Spaces: I}, Springer-Verlag, Berlin
        (1970).

\biblitem{Rudin1}
        W.~Rudin, {\it Real and Complex Analysis}, Mc Graw-Hill,  (1987).

\biblitem{Rudin3}
        W.~Rudin, {\it Functional Analysis}, Mc Graw-Hill,  (1991).

\biblitem{Brezis}
        H.~Brezis, {\it Analyse Fonctionnelle-Th\'eorie et Applications},
	Masson, Paris, (1983).

\biblitem{DB&FZ1}
        D.~Bertacchi, F.~Zucca, {\it Equidistribution of random walks
        on spheres}, J.~Stat.~Phys.~{\bf 94} (1999), 91-111.

\biblitem{Pic-Woess}     
        M.~Picardello, W.~Woess, {\it The Full Martin Boundary of The
        Bi-Tree}, Annals of Probability {\bf 22} (1994), No.~4, 2203-2222.

\biblitem{Nearest}     
        W.~Woess, {\it Nearest Neighbour Random Walks on Free Product
        of Discrete Groups}, Bollettino U.~M.~I.~(6) {\bf 5}-B (1986), 961-982.
        
\biblitem{Dixmier}     
        J.~Dixmier, {\it Les moyennes invariantes dans les
        s\'emigroupes et leurs applications} Acta Sci.~Math.~(Sze\-ged), 
        {\bf 12} A (1950), 213-227.

\biblitem{Grab-Woess}     
        P.~J.~Grabner, W.~Woess, {\it Functional iterations and
        periodic oscillations for simple random walk on the 
        Sierpi\'nski graph}, Stochastic Processes and Their Applications 
        {\bf 69} (1997), 127-138.

\biblitem{Ivanov}     
        A.~A.~Ivanov, {\it Bounding the diameter of a distance-regular 
        graph}, Soviet Math.~Dokl. {\bf 28} (1983), 149-152.

\biblitem{Survey}     
        W.~Woess, {\it Random walks on infinite graphs and groups - A 
        survey on selected topics}, Bull.~London Math.~Soc. {\bf 26} (1994),
        1-60.

\biblitem{Nechaev}     
         S.~K.~Nechaev, A.~Yu.~Grosberg, A.~M.~Vershik, {\it Random walks on
         braid groups: Brownian bridges, complexity and statistics},
         J.~Phys.~A {\bf 29} (1996), 
         no 10, 2411-2433.

\biblitem{Cassi-Regina}     
        D.~Cassi, S.~Regina, {\it Random walks on $d$-dimensional
        comb lattices}, Modern Phys.~Lett. B  {\bf 6} (1992), 1397-1403.

\biblitem{Cassi-Burioni1}     
        R.~Burioni, D.~Cassi, A.~Vezzani, {\it The type-problem 
		on the average for random walks on graphs}, preprint.

\biblitem{Avez}     
        A.~Avez, {\it Limite des quotients pour des marches al\'eatoires sur
                des groupes}, C.~R.~Acad. Sci.~Paris S\'er.~A {\bf 276}
                (1973), 317-320.

\biblitem{Avez2}     
        A.~Avez, {\it Entropie des groupes de type fini},
                C.~R.~Acad.~Sci.~Paris S\'er.~A {\bf 275}
                (1972), 1363-1366.

\biblitem{Cartwright}     
        D.~I.~Cartwright, {\it Some examples of random walks on free
                products of discrete groups}, Ann.~Mat.~Pura ed Appl.~(IV), 
                {\bf 151} (1988), 1-15.

\biblitem{Gerl}     
        P.~Gerl, {\it A local central limit theorem on some groups}, 
        in {\it The First Pannonian Symposium on Mathematical Statistics},
        Springer Lecture Notes in Statistics, {\bf 8} (1981), 73-82.

\biblitem{Pica}     
        M.~Picardello, W.~Woess, {\it Random walks on amalgams},
        Monatsh.~Math.~{\bf 100} (1985), 21-33.

\biblitem{Guivarc} Y.~Guivarc'h, {\it Sur la loi des grand nombres et le rayon spectral
        d'une marche al\'eatoire}, Ast\'erisque, {\bf 74} (1980), 47-98.

\biblitem{Woess2}
        W.~Woess, {\it Random walks on infinite graphs and groups}, 
        Cambridge Tracts in Mathematics, {\bf 138}, Cambridge Univ. Press, 2000.

\biblitem{Macpherson}
        H.~D.~Macpherson, {\it Infinite distance transitive graphs of finite
        valency}, Combinatorica, {\bf 2} (1982), 63-69.

\biblitem{Cartwright-Soardi}
        D.~I.~Cartwright, P.~M.~Soardi, 
        {\it Random walks on free products, quotients and amalgams}, 
        Nagoya Math.~j. {\bf 102} (1986), 163-180.

\biblitem{Sawyer}
        S.~Sawyer, {\it Isotropic random walks in a tree}, 
        Z.~Wahrsch.~Verw.~Gebiete, {\bf 42} (1978), 279-292.

\biblitem{Bender}
        E.~A.~Bender, {\it Asymptotic method in enumeration}, 
        SIAM review, {\bf 4} (1974), 485-515.


\def\ident{{\mathchoice {\rm 1\mskip-4mu l} {\rm 1\mskip-4mu l}
{\rm 1\mskip-4.5mu l} {\rm 1\mskip-5mu l}}}
\def\complessi{{\mathchoice {\setbox0=\hbox{$\displaystyle\rm C$}\hbox{\hbox
to0pt{\kern0.4\wd0\vrule height0.9\ht0\hss}\box0}}
{\setbox0=\hbox{$\textstyle\rm C$}\hbox{\hbox
to0pt{\kern0.4\wd0\vrule height0.9\ht0\hss}\box0}}
{\setbox0=\hbox{$\scriptstyle\rm C$}\hbox{\hbox
to0pt{\kern0.4\wd0\vrule height0.9\ht0\hss}\box0}}
{\setbox0=\hbox{$\scriptscriptstyle\rm C$}\hbox{\hbox
to0pt{\kern0.4\wd0\vrule height0.9\ht0\hss}\box0}}}}
\def\bbbe{{\mathchoice {\setbox0=\hbox{\smalletextfont e}\hbox{\raise
0.1\ht0\hbox to0pt{\kern0.4\wd0\vrule width0.3pt height0.7\ht0\hss}\box0}}
{\setbox0=\hbox{\smalletextfont e}\hbox{\raise
0.1\ht0\hbox to0pt{\kern0.4\wd0\vrule width0.3pt height0.7\ht0\hss}\box0}}
{\setbox0=\hbox{\smallescriptfont e}\hbox{\raise
0.1\ht0\hbox to0pt{\kern0.5\wd0\vrule width0.2pt height0.7\ht0\hss}\box0}}
{\setbox0=\hbox{\smallescriptscriptfont e}\hbox{\raise
0.1\ht0\hbox to0pt{\kern0.4\wd0\vrule width0.2pt height0.7\ht0\hss}\box0}}}}
\def\razionali{{\mathchoice {\setbox0=\hbox{$\displaystyle\rm Q$}\hbox{\raise
0.15\ht0\hbox to0pt{\kern0.4\wd0\vrule height0.8\ht0\hss}\box0}}
{\setbox0=\hbox{$\textstyle\rm Q$}\hbox{\raise
0.15\ht0\hbox to0pt{\kern0.4\wd0\vrule height0.8\ht0\hss}\box0}}
{\setbox0=\hbox{$\scriptstyle\rm Q$}\hbox{\raise
0.15\ht0\hbox to0pt{\kern0.4\wd0\vrule height0.7\ht0\hss}\box0}}
{\setbox0=\hbox{$\scriptscriptstyle\rm Q$}\hbox{\raise
0.15\ht0\hbox to0pt{\kern0.4\wd0\vrule height0.7\ht0\hss}\box0}}}}
\def\bbbt{{\mathchoice {\setbox0=\hbox{$\displaystyle\rm
T$}\hbox{\hbox to0pt{\kern0.3\wd0\vrule height0.9\ht0\hss}\box0}}
{\setbox0=\hbox{$\textstyle\rm T$}\hbox{\hbox
to0pt{\kern0.3\wd0\vrule height0.9\ht0\hss}\box0}}
{\setbox0=\hbox{$\scriptstyle\rm T$}\hbox{\hbox
to0pt{\kern0.3\wd0\vrule height0.9\ht0\hss}\box0}}
{\setbox0=\hbox{$\scriptscriptstyle\rm T$}\hbox{\hbox
to0pt{\kern0.3\wd0\vrule height0.9\ht0\hss}\box0}}}}
\def\bbbs{{\mathchoice
{\setbox0=\hbox{$\displaystyle     \rm S$}\hbox{\raise0.5\ht0\hbox
to0pt{\kern0.35\wd0\vrule height0.45\ht0\hss}\hbox
to0pt{\kern0.55\wd0\vrule height0.5\ht0\hss}\box0}}
{\setbox0=\hbox{$\textstyle        \rm S$}\hbox{\raise0.5\ht0\hbox
to0pt{\kern0.35\wd0\vrule height0.45\ht0\hss}\hbox
to0pt{\kern0.55\wd0\vrule height0.5\ht0\hss}\box0}}
{\setbox0=\hbox{$\scriptstyle      \rm S$}\hbox{\raise0.5\ht0\hbox
to0pt{\kern0.35\wd0\vrule height0.45\ht0\hss}\raise0.05\ht0\hbox
to0pt{\kern0.5\wd0\vrule height0.45\ht0\hss}\box0}}
{\setbox0=\hbox{$\scriptscriptstyle\rm S$}\hbox{\raise0.5\ht0\hbox
to0pt{\kern0.4\wd0\vrule height0.45\ht0\hss}\raise0.05\ht0\hbox
to0pt{\kern0.55\wd0\vrule height0.45\ht0\hss}\box0}}}}

\def\pr{{\Bbb P}}
\def\rg{{\rm Rg}}

\def\d{\,{\rm d}}

\def\phi{{\varphi}}
\def\epsilon{{\varepsilon}}
 
\def\deg{{\rm deg}}

\def\ovl{\overline}

\def\reali{{\Bbb R}}
\def\complessi{{\Bbb C}}

\def\naturali{{\Bbb N}}
\def\boxf{{$\ulcorner \! \! \lrcorner \! \! \! \! \llcorner \! \! \urcorner$}}
\def\QED{{\ifmmode\sq\else{\unskip\nobreak\hfil
\penalty50\hskip1em\null\nobreak\hfil{\boxf}
\parfillskip=0pt\finalhyphendemerits=0\endgraf}\fi} \vskip 12 pt} 

\font\fontgrande=cmss17 
\centerline{\fontgrande On some properties of transitions operators} 
\vskip 10 pt 
{\noindent 
\centerline{Fabio Zucca}\par\smallskip 
\noindent 
\centerline{Universit\`a degli Studi di Milano}\par\noindent 
\centerline{Dipartimento di Matematica F.~Enriques}\par\noindent 
\centerline{Via Saldini 50, 20133 Milano, Italy.}\par} 
\vskip 12 pt 
 
{\leftskip=80pt \rightskip=80pt \noindent {\bf Abstract.}\ 
We study a general transition operator, generated by a random walk 
on a graph $X$; in particular we give necessary and sufficient
condition on the matrix coefficient ($1$-step transition probablilities) 
to be a bounded operator from $l^\infty(X)$ into itself. Moreover
we characterize compact operators and we relate this property 
to the behaviour of the associated random walk. We give a
necessary and sufficient condition for the
pre-adjoint of the discrete Laplace operator to be an injective map.
\par 
} 
 
 
\vskip 18pt 
 
\autosez{intro}{Introduction}

In this paper we consider the transition operator $P$ 
associated with a general irreducible Markov chain $\{Z_n\}_{n\in \naturali}$ defined on a probability space $(X, \Omega, \Bbb P)$
with  a state space $X$ which is finite or countable.
It is well known that one of the consequences of the time-independence property of Markov chains is that ${\Bbb P}[Z_{n+1}=y|Z_n=x]$  is independent of $n$, hence
the transition operator $P$ is uniquely determined (according to 
equation~\eqref{operatorep}) by the coefficients $\{p(x,y)\}_{x,y\in X}$
(called {\it $1$-step transition probabilities} defined by 
$$
p(x,y):={\Bbb P}[Z_{1}=y|Z_0=x]. \autoeqno{1}
$$

This map gives rise to an amount of interesting concepts
(such as harmonic and superharmonic functions, stationary measures and so
on) which allow us to get information
about the behaviour of the Markov chains.
Take for instance the characterization of recurrent random walk in terms
of non-negative superharmonic functions (see  Theorem~1.16
of \cite{Woess2} or Theorem~5.3 of \cite{Woess3}) or in terms of
excessive measures (see \cite{Woess2} Theorem~1.18). Other 
applications may be found in the study of the asymptotic behaviour 
of the $n$-step transition probabilities 
$p^{(n)}(x,y)$ (see for instance
\cite{Coulhon-Grigoryan}).

One of the most interesting topics is given by
the discrete harmonic analysis on graphs (see \cite{Cartier2} and \cite{Fig-Picar})
and the corresponding Dirichlet problem (see for instance \cite{Woess2}, 
Chapter~4).

More recently we started to study mean value properties
for finite variation measures on graphs (see \cite{Zucca1}) with
respect to suitable families of harmonic functions;
in that paper it is shown how these properties are related
with the range of the preadjoint of the discrete Laplace operator
(see Section~\sref{kernel}). 
These are
some of the reasons which justify the present paper. For those who are
not familiar with some terminology we suggest to look at \cite{Brezis}
for functional analysis reference and at \cite{Woess2} for random walk
theory reference. 
\medskip
We begin (Section~\sref{compact}) dealing with ``generalized''
transition operator (with more general kernels $p(x,y)$) defined
by the equation~\eqref{operatorep}: we give a complete 
characterization of continuous transition maps
(Theorem~\lemmaref{boundP}) and of compact transition
maps with non-negative kernel (Theorem~\lemmaref{compact1});
moreover we show that compactness implies recurrence 
(Proposition~\lemmaref{recurrentcom}), meanwhile
the converse it is not true (Example~\lemmaref{es1}).

In Section~\sref{kernel} we consider the stochastic kernel defined by
equation~\eqref{1} and we deal with the corresponding
(stochastic) transition operator
which is a continuous map from $l^\infty(X)$ into itself. We construct the 
preadjoint of the discrete analogous of the Laplace operator and 
we turn our attention to its null space and its range. In particular
we give a necessary and sufficient condition
for this operator to be an injective
map (Theorem~\lemmaref{stationary}); this result generalizes 
Theorem~1.18 of \cite{Woess2}.
We finally make some remarks about the topological properties of its
range.
\medskip

We fix now the basic notation:
let $\Phi$ represent the real field  $\reali$ or the complex
field $\complessi$ and
$p:X \times X \rightarrow \Phi$ be a function. 
We consider the domain $D$ and the linear operator $P$
depending by $p$  as follows
$$
\eqalign{
D(P)&:=\left\{f:X \rightarrow \Phi : 
\sum_{y \in X} |p(x,y)| |f(y)| <+\infty, \ \forall x\in X \ \right\}, \cr
(Pf)(x)&:=\sum_{y \in X} p(x,y) f(y), \qquad \forall f \in D(P), \ \forall
x \in X.\cr
} \autoeqno{operatorep}
$$

The properties of the linear map
$P$ are strictly related to
the functional space where it is restricted:
for instance if the coefficients $p$ satisfy equation~\eqref{1}
(which is equivalent to
$p(x,y)\geq 0$ for all $x,y \in X$ and $\sum_{y \in X} p(x,y)=1$ for
every $x\in X$) then the transition operator
$P$ is called {\it stochastic}; in this case it is easy to show that
$P$ is a bounded linear map from $l^\infty(X, \mu)$ into
itself (for any real or complex measure $\mu$ on $X$) and $\|P\|=1$;
furthermore given any excessive, positive measure $\nu$ on $X$ (see
Section~\sref{kernel}),
any stochastic map $P$ is bounded, with $\|P\| \leq 1$, from
$L^p(X, \nu)$ into itself ($p \in [1, +\infty)$).

If $P$ is generated by a reversible random walk
(see for instance \cite{Survey}, Paragraph 2.A) and if $\nu$ is the reversibility measure then
$P$ is a linear, bounded, selfadjoint operator from
$L^2(X,\nu)$ into itself, satisfying $\|P\|_{L^2(X,\nu)} = \rho(P)$,
where $\rho(P)$ is the spectral radius of the random walk $(X,P)$ (see \cite{Woess2} Chapter 1, Paragraph B).
More precisely it is possible to show that an operator $K$ defined as in
eq.~\eqref{operatorep} with kernel $k(x,y)$ (instead of $p(x,y)$)
is selfadjoint if and only if
$\nu(x)k(x,y)=\nu(y)k(y,x)$ for all $x,y \in X$.

\autosez{compact}{Compactness of the transition operator}

In this section we give a necessary and sufficient condition
for the general linear map $P$ (defined by eq.~\eqref{operatorep}) to be
a bounded map from 
from $l^\infty(X)$ into itself.
Then we characterize all the maps with non-negative kernels which are
compact; in the case of stochastic maps, 
this conditions is related to
the recurrence property. The interest in the space $l^\infty(X)$ will be 
justified in the next section.
We start
with boundeness conditions.

\theorem{boundP}{Let $P$ the transition operator defined by
eq.~\eqref{operatorep} (where $p(x,y)$ are real (complex) numbers for
any $x,y \in X$); then the following assertions are equivalent:
\vskip 6 pt
\item{(i)} $P$ is a continuous linear operator from $l^\infty(X)$
into itself; 
\vskip 6 pt
\item{(ii)} $\sup_{x \in X} \sum_{y \in X} |p(x,y)| < \infty$. 
\vskip 6 pt
\item{(iii)} $D(P) \supseteq l^\infty(X)$ and $P(l^\infty(X)) \subseteq 
l^\infty(X)$.
\vskip 6 pt
\noindent If one of the previous condition holds, then
$\|P\|= \sup_{x \in X} \sum_{y \in X} |p(x,y)|$.
}
\proof Let us discuss the complex case.
Since $X$ is at most countable, we may suppose $X = \{x_i\}_{i \in J}$, where $J \subseteq \naturali$ has the same cardinality of $X$ . If
$X$ is finite then the theorem is trivial, we suppose that $X$ is countable (and $J=\naturali$).
Let $x \in X$ and define
$$
f^n_x(x_i):= \cases{ {\overline{p(x,x_i)} \over |p(x,x_i)|} & if $p(x,x_i) \not = 0$,
$i\leq n$, \cr
\cr
1 & if $p(x,x_i)=0$, $i \leq n$, \cr
\cr
0 & if $i > n$. \cr
}
$$
then $\|f^n_x\|_\infty =1$ for every $n \in \naturali$ 
and for every $x \in X$.
\smallskip
(i) $\Longrightarrow$ (ii). We easily note that,
 $f_x^n \in D(P)$ for every $n\in\naturali, x \in X$ and
$\|Pf_x^n\|_\infty \geq |(Pf_x^n)(x)|= \sum_{i \leq n} |p(x,x_i)|$,
then if $n$ tends to infinity and $P$ is bounded, we have
$$
\sum_{y \in X} |p(x,y)| = \sum_{i \in \naturali} |p(x,x_i)|  \leq \|P\|
\autoeqno{compact1}
$$
and the arbritrary choice of $x$ leads to the conclusion.

(ii) $\Longrightarrow$ (i).
If $\alpha := \sup_{x \in X} \sum_{y \in X} |p(x,y)| < \infty$
then for every $f \in l^\infty(X)$ we obtain
$$
\|Pf\|_\infty = \sup_{x \in X} |(Pf)(x)| \leq
\sup_{x \in X} \sum_{y \in X} |p(x,y)| \, \|f\|_\infty= \alpha \|f\|_\infty;
$$
using the last equation and  eq.~\eqref{compact1} we obtain
that $D(P) \supseteq l^\infty(X)$, $P(l^\infty(X)) \subseteq l^\infty(X)$
and $\alpha = \|P\|$.

(iii) $\Longrightarrow$ (ii). 
Let us note that $l^\infty(X) \subseteq D(P)$ implies
that for every $f \in l^\infty(X)$ and for every $x \in X$ we have
$\sum_{y \in X}|p(x,y)||f(y)| < +\infty$ which is equivalent to
the condition $\{p(x, y)\}_{y\in X} \in l^1(X)$ for every $x \in X$.
Let $\lambda_x \in l^\infty(X)^*$ defined by $\lambda_x(f):=\sum_{y \in X}
p(x,y)f(y)$, then $\|\lambda_x\|_{l^\infty(X)^*} = \sum_{y \in X} |p(x,y)|$.
Now the condition $P(l^\infty(X)) \subseteq l^\infty(X)$ implies
$\sup_{x \in X} |\lambda_x(f)| < +\infty$, for every $f \in l^\infty(X)$, then,
using the principle of Uniform Boundeness, we have 
$\sup_{x \in X}  \|\lambda_x\|_{l^\infty(X)^*} < +\infty$ which is
equivalent to (ii).

(i) $\Longrightarrow$ (iii). It is trivial.

\QED

We turned our attention now to the compactness property for a transition operator
with non negative kernel.

\theorem{compact1}{Let $X$ a countable graph and
let us choose an enumeration $\{x_i\}$ for $X$. Let $P$ a transition operator
on $X$ with non negative elements, satisfying the condition
$\sup_{x \in X} \sum_{y \in X} p(x,y) <+\infty$. Then $P$ is a bounded,
linear operator from $l^\infty(X)$ into itself; moreover $P$ is compact
if and only if
$$ \lim_{n \rightarrow +\infty} \sup_{x \in X} \sum_{i> n} p(x,x_i) =0.
\autoeqno{compact0}
$$
The last condition is independent from the chosen enumeration.
}
\proof
Theorem~\lemmaref{boundP} implies the boundedness of $P$. 
If $P$ is compact, then, for every bounded sequence $\{f_i\}$ in
$l^\infty(X)$, $\{Pf_i\}$ is relatively compact, hence
there exists a subsequence $\{n_j\}$ such that $\{Pf_{n_j}\}$ is a
Cauchy sequence. Let $f_n(x_i)$ equal to $1$ if $i >n$ and $0$ otherwise,
then
$$
(Pf_n)(\cdot)= \sum_{i >n } p(\cdot, x_i)
$$ and if $m >n$, since $p(x,y) \geq 0$ for every $x,y \in X$,
$$
\|Pf_n-Pf_m\|_\infty = \sup_{x \in X} \sum_{i=n+1}^m p(x,x_i).
$$
By the Cauchy property, for every $\epsilon >0$ there exists $j_\epsilon$
such that for all $j_2 > j_1 \geq j_\epsilon$ we have
$$
\sup_{x \in X} \sum_{i=n_{j_1}+1}^{n_{j_2}} p(x,x_i) < \epsilon/2;
$$ we note
that, for every fixed $x \in X$, $m \mapsto \sum_{i =n+1}^m p(x,x_i)$
(resp.~$n \mapsto \sum_{i =n+1}^m p(x,x_i)$) is not decreasing
(resp.~not increasing), hence
$$
\sup_{x \in X} \sum_{i=n_{j_\epsilon}+1}^\infty p(x,x_i) \leq \epsilon/2<
\epsilon,
$$
which implies $\lim_{n \rightarrow +\infty} \sup_{x \in X} \sum_{i> n} p(x,x_i) =0$.

Vice versa if we consider the finite range (compact) projections
on $l^\infty(X)$ defined by
$$
V_i(f)(x_n):=\cases{ f(x_n) & if $n \leq i$ \cr \cr
	0 & if $n > i$, \cr
}$$
for all $i \in \naturali$, then
$(PV_if)(x)= \sum_{i \leq n} p(x,x_i) f(x_i)$.
By Theorem~\lemmaref{boundP} 
$\|P-PV_i\|= \sup_{x \in X} \sum_{i >n} p(x,x_i)$ then
$\|P-PV_i\|$ tends to $0$ if $n$ tends to infinity. By
Theorem~4.18(f) of Rudin~\cite{Rudin1}, $PV_i$ is compact for every
$i \in \naturali$, hence $P$ is compact since
Theorem~4.18(c) of Rudin~\cite{Rudin1} holds.

The condition~\eqref{compact0} does not depend on the choice
of the enumeration, since compactness is defined ``a priori''.
\QED

We note that the previous theorem says that $P$ is compact
if and only if for any
$\epsilon >0$ there exists a finite subset $A_\epsilon \subset X$
such that  $\sup_{x \in X} \sum_{y\in A_\epsilon} p(x,y) < \epsilon$;
this means, for instance, that a necessary condition for the
compactness property is that $\lim_{y \rightarrow \infty} p(x,y)=0$
holds uniformly with respect to $x \in X$ (where the limit is taken in
the Alexandroff compactification of $X$ with the discrete topology).

We immediately note that if $X$ is locally finite, then $P$ is not compact.
In fact in this case, for every $n \in \naturali$, $\sum_{i \leq n}
\deg(x_i) <+\infty$; this means that there exists $n_1 >n$ such that 
$x_i$ is not a neighbour of $x_{n_1}$ for any $i\leq n$, hence
$\sum_{i>n} p(x_{n_1},x_i) =1$ and eq.~\eqref{compact0} cannot be satisfied.
Moreover one can show that if $P$ is compact then it is a recurrent 
transition operator. 

\proposition{recurrentcom}{Let $X$ be a graph and $P$ a random walk on $X$
which is a compact,
transition operator from $l^\infty(X)$ into itself; then $P$ is recurrent.
}
\proof
Let $A:=\{x_0,x_1, \ldots, x_n\}$ and $Z_n$ the Markov chain associated to 
$P$; when $P$ is compact then, by Theorem~~\lemmaref{compact1}, if $x \in X$
$$
\pr [Z_m \not \in A | Z_0 =x]= \sum_{\scriptstyle y \in X \atop i > n}
\pr [Z_{m-1}=y|z_0=x] p(y,x_i)={}
$$
$$
{}=\sum_{y \in X} \pr [Z_{m-1}=y|z_0=x] 
\sum_{i > n} p(y,x_i)\leq {}
$$
$$
{}\leq
\sum_{y \in X} \pr [Z_{m-1}=y|z_0=x] 
\sup_{w \in X}\sum_{i > n} p(w,x_i)\leq
\sup_{w \in X} \sum_{i>n} p(w,x_i) \ 
{\buildrel n \rightarrow +\infty \over \longrightarrow} \ 0
$$
If $B:=\{\exists k \in \naturali; Z_n \not \in A, \, \forall n \geq k\}
\equiv \cup_{k \in \naturali} \cap_{n \geq k} \{Z_n \not \in A\}$
it is clear that 
$$\pr (B) \leq \sum_{k \in \naturali} \pr (\cap_{n \geq k} \{Z_n \not \in A\}),
$$
but $\cap_{n \geq k} \{Z_n \not \in A\} \subseteq \{Z_m \not \in A\}$
for every $m \geq k$ which implies 
$\pr (\cap_{n \geq k} \{Z_n \not \in A\})=0$ and $\pr (B)=0$. Thus by Woess~\cite{Woess3}
Theorem~3.5, $P$ is recurrent.
\QED

The necessary condition given in the previous proposition is not sufficient; it
is not difficult to find random walks which are recurrent and the associated transition
operator is not compact.

\example{es1}{Let us give an example of random walk giving rise to a
compact transition operator. To this aim we consider any sequence 
of real number $\{p_i\}_{i\in \naturali}$ such that $p_0=1$ and 
$p_i\in
(0,1]$ for
every $i \geq 1$. Let us take $X=\naturali$ and
$$
p(x,y):=\cases{
p_x & if $x \in \naturali$ and $y=x+1$ \cr
1-p_x & if $y=0$ and $x\not = 0$ \cr
0 & otherwise. \cr
}
$$
The condition~\eqref{compact0} is easily equivalent to 
$\lim_{n \rightarrow +\infty} p_n =0$; moreover, using
Theorem 3.11 of \cite{Woess3}, it is not
difficult to show that $(X,P)$ is recurrent (resp.~positive
recurrent) if and only if $\lim_{n \rightarrow +\infty} \prod_{i=0}^n
p_i =0$ (resp.~$\sum_{n=0}^\infty \prod_{i=0}^n p_i < +\infty$).
This proves that $(X,P)$ positive recurrent (and hence $(X,P)$ recurrent)
does not imply the compactness of the transition operator $P$.
}

\autosez{kernel}{The null space of the pre-adjoint of the Laplace operator and
finite variation stationary measures.}

The discrete analogous of the Laplace operator is given by
$(P-\ident_\infty):l^\infty(X) \rightarrow l^\infty(X)$
where $\ident_\infty$ id the identity operator on $l^\infty(X)$. It is
straightforward to show that the preadjoint map
$(P-\ident_\infty)_*$ is given by 
$(Q-\ident_1):l^1(X) \rightarrow l^1(X)$ where
$$
(Q\nu)(y) :=\sum_{x \in X} p(x,y) \nu(x), \qquad
\forall y \in X,
$$
and $\ident_1$ is the identity map on $l^1(X)$; from now on, 
if $A:D \rightarrow Y$ is a map, we denote by $\rg(A)$ the image (or range)
of $D$ through $A$ (that is $A(D)$). A bounded function is said to be {\it
harmonic} if is an element of the null space of the discrete laplacian
(we denote the set of all bounded harmonic functions by ${\cal 
H}^\infty(X,P)$).

In \cite{Zucca1} it was shown that a finite variation measure $\nu$ on $X$ (i.e.
~$\nu \in l^1(X)$) has the weak mean value property with respect to $o \in X$
(that is, $\sum_{x \in X} f(x) \nu(X) \d \nu =f(o) \sum_{x \in X} \nu(x)$,
for every $f \in {\cal H}^\infty(X,P)$) if and only if
$(\nu-\delta_o \sum_{x \in X} \nu(x)) \in \overline{\rg(Q-\ident_1)}$
(where $\delta_o$ is the (finite variation) Dirac measure with support in $\{o\}$.
>From this point of view, it is important to know when $Q-\ident_1$ is injective 
and when it has a closed range. In this section we give a complete
answer to the first question and we make some remarks related to the second.

\smallskip
To this aim, we characterize in particular all the stationary measures
with finite variation. We recall that a signed measure is called
{\it stationary} (resp.~{\it excessive}) if
$$
(Q\nu)(y)=\nu(y), \quad \forall y \in X
 \qquad \hbox{(resp.~}(Q\nu)(y) \leq \nu(y), \quad \forall y \in X \hbox{)}
$$
provided that $(Q \nu)(y)$ exists for every $y \in X$.
We note that a finite variation measure $\nu$ is stationary
if and only if $\nu \in \ker(Q-\ident_1)$.

\lemma{totvar}{Let $\nu$ be a signed stationary measure on $X$, then
$-|\nu|$ is a negative excessive measure which is stationary if and only if
$\nu = |\nu|$ or $\nu= -|\nu|$.}
\proof
It is well know that if $f$ is a complex integrable function on
a measure space $(Y, \mu)$ then
$$
\left | \int_Y f \d \nu \right | \leq \int_Y |f| \d \nu
$$
and the equality holds if and only if there exists $\alpha \in [0, 2 \pi )$
such that $f = |f| \exp(i \alpha)$ $\mu$-a.e.
If we consider the measure space $X$ with the counting measure and
$f_y(x):=\nu(x)p(x,y)$ then by hypothesis $f_y \in L^1(X)$ for every
$y \in X$ and
$$
|\nu|(y)=|\nu(y)| = \left | \sum_{x \in X} \nu(x) p(x,y) \right | \leq
\sum_{x \in X} |\nu|(x) p(x,y) = (Q|\nu|)(y)
$$
and the equality holds if and only if
$\nu(x) p(x,y) = |\nu|(x) p(x,y) \exp(i \alpha)$ (where $\alpha \in
\{o, \pi \}$, since $f_y$ is a real function) which leads to the conclusion.
\QED

We are ready to prove the following theorem which characterizes
all the stationary measures (i.e.~the null space of $Q-{\Bbb I}_1$).
Before state the Theorem we recall that 
for any irreducible random walk it is possible to define
a natural number called the {\it period} of the random walk
(see \cite{Survey}, Section~5.A).

\theorem{stationary}{Let $(X,P)$ be an irreducible random walk, then there
exists a finite variation, stationary measure $\nu \not \equiv 0$ if
and only if $(X,P)$ is positive recurrent. In this case there exists
$\alpha \in \reali \setminus \{o\}$ such that $\nu = \alpha \mu$
where $\mu$ satisfies 
$$ \mu(y)= \limsup_{n \rightarrow \infty} p^{(nd +j-i)}(x,y)/d
$$
(the right hand side is seen to be independent of $x$ and $d$ is 
the period of the random walk).
}
\proof
If we suppose that there exists a stationary measure $\nu$ with
finite variation and $C_0, C_1, \dots , C_{d-1}, C_d \equiv C_0$
is the partition of $X$ given by the
periodicity classes, then, by 
Lemma~\lemmaref{totvar}, $-|\nu|$
is an excessive measure; 
Lemmas~2.4, 2.5 and Theorem~2.2 of Woess~\cite{Woess3}
 and Tonelli-Fubini's
Theorem imply
$$
|\nu|(C_{i+1}) \sum_{y \in C_{i+1}} |\nu(y)|
\leq \sum_{y \in C_{i+1}} \sum_{x \in C_i} |\nu(x)| p(x,y) =
{}
$$
$$
{} = \sum_{x \in C_i} \sum_{y \in C_{i+1}} |\nu(x)| p(x,y) =
|\nu|(C_i)
$$
then $|\nu|(C_0) \leq |\nu|(C_1) \leq \cdots \leq |\nu|(C_{d-1})
\leq |\nu|(C_d) \equiv |\nu|(C_0)$, hence
$|\nu|(C_i) = |\nu|(X) / d$ for every $i =0,1, \ldots, d-1$.

Using the ``Renewal Theorem'' by Erd\"os-Feller-Pollard 
(see \cite{ErFePo} or Theorems~3.6 and 3.7 of \cite{Woess3})
 and Lebesgue bounded convergence Theorem, 
$$
|\nu|(y) \leq \sum_{x \in C_i} |\nu|(x) p^{(nd)}(x,y) \,
{\buildrel n \rightarrow +\infty \over \longrightarrow}
d \cdot \mu(y) \sum_{x \in C_i} |\nu|(x); \autoeqno{dis1}
$$
since $\nu \not \equiv 0$ then there exists $i$ and $y \in C_i$
such that $|\nu|(y) > 0$, thus eq.~\eqref{dis1}
implies that $\mu(y) >0$, hence $(X,P)$ is positive recurrent.

On the other hand, if $(X,P)$ is positive recurrent, Theorem~3.9
of Woess~\cite{Woess3} implies that $\mu$ is a stationary, probability measure.

If $\nu$ is another stationary, finite variation measure on $X$
($\nu\not \equiv 0$) then by eq.~\eqref{dis1},
$|\nu|(y) \leq |\nu|(X) \mu(y)$; if we suppose, by contradiction, that
there exist $y \in X$ such that
$|\nu|(y) < |\nu|(X) \mu(y)$, then we have
$$
1 = \sum_{y \in X} |\nu|(y) /|\nu|(X) < \sum_{y \in X} \mu(y)=1
$$
which is a contradiction; hence $|\nu|(\cdot) / |\nu|(X) \equiv \mu(\cdot)$.
If we define $\ovl \nu (y) := \nu(y) / |\nu|(X)$ then
$|\ovl \nu | \equiv \mu$ and
$(2 \mu - \ovl \nu) / (2 \mu(X) - \ovl \nu(X))$
 is a stationary, probability measure.
By Theorem~3.9 of Woess~\cite{Woess3}, 
$(2 \mu - \ovl \nu) / (2 \mu(X) - \ovl \nu(X)) \equiv \mu$
which is equivalent to $\ovl \nu= \ovl \nu(X) \mu$, that is,
$\nu = \nu(X) \mu$.
\QED

As a consequence of this Theorem we obtain that the bounded, linear map
$Q- \ident_1$ is injective if and only if $(X,P)$ is
not positive recurrent.

We try now to answer to the second question: when $\rg(Q-\ident_1)$ is closed?

By Schauder's Theorem (see Brezis~\cite{Brezis} Theorem~VI.4), since $P=Q^*$,
 we have that the operator $Q$ from $l^1(X)$
into itself, is compact if and only if
eq.~\eqref{compact0} holds. Now it is well known 
(see for instance \cite{Rudin3}, Theorem~4.23) that if $Q$ is
compact operator from a Banach space into itself then
$Q-\ident$ has closed range. Hence if equation~\eqref{compact0}
holds we have
that a finite variation measure $\nu$ on $X$ 
has the weak mean value property with respect to $o \in X$
if and only if
$(\nu-\delta_o \sum_{x \in X} \nu(x)) \in \rg(Q-\ident_1)$.

This is obviously 
only a partial answer to our question; a way
to reach a complete a satisfactory answer,
which we don't undertake here is given by
the following remarks.

We recall that if $(Z, \|\cdot\|_Z)$, $(Y,\|\cdot\|_Y)$ are Banach spaces,
$D$ is a linear subspace of $Z$ and $A:D \rightarrow Y$
is a linear map such that $\sup_{x \in D: \|x\|_Z=1} \|Ax\|_Y=:\beta$,
then there exists a unique bounded, linear map $\overline{A}: \overline{D} \rightarrow Y$
which extends $A$; moreover $\overline{A}$ is bounded by the same constant $\beta$
and
$\displaystyle \inf_{x\in D: \|x\|_Z=1} \|Ax\|_Y = \inf_{x\in \overline{D}: \|x\|_Z=1}
 \|\overline{A}x\|_Y$.
Therefore, if $A:D \rightarrow Y$ is a linear and injective map, then
$$
\sup_{y \in \rg(A): \|y\|_Y=1} \|A^{-1}y\|_Z
 = 1/ \inf_{x\in D: \|x\|_Z=1} \|Ax\|_Y
$$
(where, by definition, $1/0:=+\infty$).

Now using the Open Mapping Theorem it is simple to show that 
if $A:Z \rightarrow Y$ is a linear, bounded, injective map then
$$
\rg(A)=\overline{\rg(A)} \Longleftrightarrow 
\inf_{x\in D: \|x\|_Z=1} \|Ax\|_Y >0.
$$

In our case, if the Markov chain is not positive recurrent, 
then $\rg(Q-\ident_1)$ is closed if and only if 
$\displaystyle \inf_{\nu \in l^1(X): \|\nu\|_1=1} \|Q\nu-\nu\|_1 >0$.

\vskip 40 pt 
{\bf Bibliography} 
\bigskip 
\insertbibliografia

\end